\documentclass[11pt]{article}
\usepackage{amssymb}
\usepackage{amsfonts}
\usepackage{graphicx}
\usepackage{amsmath}
\usepackage[francais]{babel}

\newtheorem{theorem}{Theorem}

\newtheorem{proposition}[theorem]{Proposition}

\newtheorem{theoreme}[theorem]{Th\'eor\`eme}

\newtheorem{corollaire}[theorem]{Corollaire}
\newtheorem{exemple}[theorem]{Exemples}
\newtheorem{remarque}[theorem]{Remarque}

\newenvironment{demonstration}[1][D\'emonstration]{\textbf{#1:} }
{\ \rule{0.5em}{0.5em}}
\input{tcilatex}

\begin{document}

\title{Parall\'{e}lisme d'une vari\'{e}t\'{e} des points proches}
\author{Basile Guy Richard BOSSOTO \\
Universit\'{e} Marien NGOUABI, Facult\'{e} des Sciences,\\
D\'{e}partement de Math\'{e}matiques\\
B.P.69, Brazzaville, Congo.\\
E-mail: bossotob@yahoo.fr}
\date{}
\maketitle

\begin{abstract}
On consid\`{e}re $M$ une vari\'{e}t\'{e} diff\'{e}rentielle, $A$ une alg\`{e}%
bre locale, $M^{A}$ la vari\'{e}t\'{e} des points proches de $M$ d'esp\`{e}%
ce $A$. On utilise la structure de $C^{\infty }(M^{A},A)$-module sur
l'ensemble $\mathfrak{X}(M^{A})$ des champs de vecteurs sur $M^{A}$ pour
donner l'\'{e}quivalence du parall\'{e}lisme de $M^{A}$ en termes de $A$-vari%
\'{e}t\'{e}s.

\textbf{Summary: }Let $M$ be a smooth manifold, $A$ a local algebra, $M^{A}$
the manifold of near points on $M$ of kind $A$. We use the structure of $%
C^{\infty }(M^{A},A)$-module on the set $\mathfrak{X}(M^{A})$ of vector
fields on $M^{A}$ \ for to give the equivalence of parallelism of the $A$%
-manifold $M^{A}$.
\end{abstract}

\textbf{Key words}: Point proche, alg\`{e}bre locale, vari\'{e}t\'{e} parall%
\'{e}lisable.

\textbf{MSC (2000)}: 58A20, 58A32.

\section{ Introduction}

Une alg\`{e}bre locale (au sens de Weil) est une alg\`{e}bre r\'{e}elle
commutative unitaire de dimension finie ayant un id\'{e}al maximal unique de
codimension $1$ sur $%
\mathbb{R}
$.

Soit $A$ une alg\`{e}bre locale et soit $\mathfrak{m}$ son unique id\'{e}al
maximal. On a 
\begin{equation*}
A=%
\mathbb{R}
\oplus \mathfrak{m}\text{.}
\end{equation*}

La premi\`{e}re projection 
\begin{equation*}
A=%
\mathbb{R}
\oplus \mathfrak{m}\longrightarrow 
\mathbb{R}%
\end{equation*}%
est un homomorphisme d'alg\`{e}bres qui est surjectif, appel\'{e} augmention
et l'unique entier naturel $k\in 
\mathbb{N}
$ tel que $\mathfrak{m}^{k}\neq (0)$ et $\mathfrak{m}^{k+1}=(0)$ est la
hauteur de $A$.

Comme exemples d'alg\`{e}bres locales, on a :

\begin{exemple}
\begin{description}
\item[1-] $%
\mathbb{R}
=%
\mathbb{R}
\oplus (0)$ est une alg\`{e}bre locale de hauteur $0$.

\item[2-] L'alg\`{e}bre des nombres duaux, $\mathbb{D}=%
\mathbb{R}
\left[ T\right] /(T^{2})$, est une alg\`{e}bre locale de hauteur $1$.

\item[3-] $\mathbb{A}=%
\mathbb{R}
\left[ T\right] /(T^{3})$ est une alg\`{e}bre locale de hauteur $2$. Plus g%
\'{e}n\'{e}ralement, l'alg\`{e}bre des polyn\^{o}mes tronqu\'{e}e%
\begin{equation*}
A=%
\mathbb{R}
\left[ X_{1},...,X_{n}\right] /(X_{1},...,X_{n})^{k+1}
\end{equation*}%
est une alg\`{e}bre locale de hauteur $k$.

\item[4-] Si $A$ est une alg\`{e}bre locale d'id\'{e}al maximal $\mathfrak{m}%
_{A}$ de hauteur $h$ et si $B$ est une alg\`{e}bre locale d'id\'{e}al
maximal $\mathfrak{m}_{B}$ de hauteur $l$, alors le produit tensoriel $%
A\otimes B$ est une alg\`{e}bre locale d'id\'{e}al maximal $\mathfrak{m}%
_{A\otimes B}=\mathfrak{m}_{A}\otimes B+A\otimes \mathfrak{m}_{B}$ et de
hauteur $h+l$. Ainsi, $\mathbb{D\otimes D=}%
\mathbb{R}
\left[ T_{1},T_{2}\right] /(T_{1}^{2},T_{2}^{2})$ est une alg\`{e}bre locale
de hauteur $2$.
\end{description}

\begin{remarque}
Le produit tensoriel de deux alg\`{e}bres de polyn\^{o}mes tronqu\'{e}es
n'est pas une alg\`{e}bre de polyn\^{o}mes tronqu\'{e}e. Ce qui est le cas
pour $\mathbb{D\otimes D}$.
\end{remarque}

\begin{description}
\item[5-] Si $M$ est une vari\'{e}t\'{e} diff\'{e}rentielle de dimension $n$%
, l'espace, $J_{x}^{k}(M,%
\mathbb{R}
)$, des jets en $x\in M$ d'ordre $k$ des applications diff\'{e}rentiables de
classe $C^{\infty }$ d\'{e}finies au voisinage de $x$ \`{a} valeurs dans $%
\mathbb{R}
,$ est une alg\`{e}bre locale de de dimension $\complement _{n+k}^{k}$ et de
hauteur $k$.
\end{description}
\end{exemple}

\bigskip

Si $M$ est une vari\'{e}t\'{e} diff\'{e}rentielle, $C^{\infty }(M)$ l'alg%
\`{e}bre des fonctions num\'{e}riques sur $M$ et $A$ une alg\`{e}bre locale,
un point proche de $x\in M$ d'esp\`{e}ce $A$ est un homomorphisme d'alg\`{e}%
bres%
\begin{equation*}
\xi :C^{\infty }(M)\ \ \longrightarrow A
\end{equation*}%
tel que $[\xi (f)-f(x)]\in m$ pour tout $f\in C^{\infty }(M)$.

On note $M_{x}^{A}$ l'ensemble des points proches de $x$ d'esp\`{e}ce $A$ et 
\begin{equation*}
M^{A}=\bigcup\limits_{x\in M}M_{x}^{A}\text{.}
\end{equation*}

L'ensemble $M^{A}$ est une vari\'{e}t\'{e} diff\'{e}rentielle de dimension $%
\dim M\cdot \dim A$.

\begin{exemple}

\begin{description}
\item[1-] $M^{%
\mathbb{R}
}=M$.

\item[2-] Pour toute vari\'{e}t\'{e} diff\'{e}rentielle $M$, l'application 
\begin{equation*}
TM\longrightarrow Hom_{A\lg }(C^{\infty }(M),\mathbb{D}),v\longmapsto \xi
_{v},
\end{equation*}%
\qquad\ d\'{e}finie par 
\begin{equation*}
\xi _{v}(f)=f(p)+v(f)\varepsilon
\end{equation*}%
si $v\in T_{p}M$, identifie $TM=J_{0}^{1}(%
\mathbb{R}
,M)$ \`{a} $M^{\mathbb{D}}.$ On v\'{e}rifie que $v$ est un vecteur tangent
en $p$ \`{a} la vari\'{e}t\'{e} $M$ si et seulement si $\xi _{v}$ est un
point proche de $p$ d'esp\`{e}ce $\mathbb{D}$.

\item[3-] Si $A=%
\mathbb{R}
\lbrack X]/(X^{3})$, $M^{A}=J_{0}^{2}(%
\mathbb{R}
,M)$. Plus g\'{e}n\'{e}ralement, si $A$ est l'alg\`{e}bre des polyn\^{o}mes
tronqu\'{e}s 
\begin{equation*}
\ 
\mathbb{R}
\lbrack X_{1},...,X_{s}](X_{1},...,X_{s})^{k+1}\text{,}
\end{equation*}
alors $M^{A}=J_{0}^{k}(%
\mathbb{R}
^{s},M)$ est l'ensemble des jets en $0$ d'ordre $k$ des applications diff%
\'{e}rentiables de $%
\mathbb{R}
^{s}\ $dans $M$.

\item[4-] L'application $\xi \longmapsto \xi (id_{%
\mathbb{R}
})$ identifie $%
\mathbb{R}
^{A}$ \`{a} $A$.

\item[5-] Si $V$ est un espace vectoriel r\'{e}el de dimension finie, si $%
(e_{i})_{i=1,...,r}$ est une base de $V$ et si $(e_{i}^{\ast })_{i=1,...,r}$
est la base duale de la base $(e_{i})_{i=1,...,r}$, alors 
\begin{equation*}
V^{A}\longrightarrow V\otimes A,\xi \longmapsto \overset{r}{\underset{i=1}{%
\dsum }}e_{i}\otimes \xi (e_{i}^{\ast })
\end{equation*}%
est un isomorphisme canonique de $A$-modules.
\end{description}
\end{exemple}

Lorsque $M$ et $N$ sont deux vari\'{e}t\'{e}s diff\'{e}rentiables et lorsque 
\begin{equation*}
h:M\longrightarrow N
\end{equation*}%
est une application diff\'{e}rentiable de classe $C^{\infty }$, alors
l'application 
\begin{equation*}
h^{A}:M^{A}\longrightarrow N^{A},\xi \longmapsto h^{A}(\xi ),
\end{equation*}%
telle que, pour tout $g\in C^{\infty }(N)$, 
\begin{equation*}
\left[ h^{A}(\xi )\right] (g)=\xi (g\circ h)
\end{equation*}%
est diff\'{e}rentiable de classe $C^{\infty }$. Lorsque $h$ est un diff\'{e}%
omorphisme, il en est de m\^{e}me de $h^{A}$.

De plus, si $\varphi :A\longrightarrow B$ est un homomorphisme d'alg\`{e}%
bres locales, pour toute vari\'{e}t\'{e} diff\'{e}rentielle $M$,
l'application 
\begin{equation*}
\varphi _{M}:M^{A}\longrightarrow M^{B},\xi \longmapsto \varphi \circ \xi
\end{equation*}%
est diff\'{e}rentiable. En particulier, l'augmention 
\begin{equation*}
A\longrightarrow 
\mathbb{R}%
\end{equation*}%
d\'{e}finit pour toute vari\'{e}t\'{e} $M$, la projection 
\begin{equation*}
M^{A}\longrightarrow M,
\end{equation*}%
qui a un point proche de $x$ $\in M$, associe son origine $x$.

\section{Parall\'{e}lisme de la vari\'{e}t\'{e} des points proches}

\bigskip Dans tout ce qui suit $M$ d\'{e}signe une vari\'{e}t\'{e} diff\'{e}%
rentielle de dimension $n$, $A$ une alg\`{e}bre locale au sens de Weil, d'%
\'{e}l\'{e}ment unit\'{e} $1_{A}$, $C^{\infty }(M)$ l'alg\`{e}bre des
fonctions num\'{e}riques de classe $C^{\infty }$ sur $M$, $\mathfrak{X}(M)$
le $C^{\infty }(M)$-module des champs de vecteurs sur $M$, $TM$ le fibr\'{e}
tangent \`{a} $M$ et 
\begin{equation*}
\pi _{M}:TM\longrightarrow M
\end{equation*}%
la projection canonique.

Si $(U,\varphi )$ est une carte locale de $M$ de fonctions coordonn\'{e}es $%
(x_{1},x_{2},...,x_{n})$, l'application,%
\begin{equation*}
U^{A}\longrightarrow A^{n},\xi \longmapsto (\xi (x_{1}),\xi (x_{2}),...,\xi
(x_{n})),
\end{equation*}%
est une bijection de $U^{A}$ sur un ouvert de $A^{n}$. La vari\'{e}t\'{e} $%
M^{A}$ est une vari\'{e}t\'{e} model\'{e}e sur $A^{n}$, c'est-\`{a}-dire une 
$A$-vari\'{e}t\'{e} de dimension $n$.

L'ensemble, $C^{\infty }(M^{A},A)$, des fonctions de classe $C^{\infty }$
sur $M^{A}$ \`{a} valeurs dans $A$, est une $A$-alg\`{e}bre commutative
unitaire. En identifiant $%
\mathbb{R}
^{A}$ \`{a} $A$, pour $f\in C^{\infty }(M)$, l'application 
\begin{equation*}
f^{A}:M^{A}\longrightarrow A,\xi \longmapsto \xi (f)\text{,}
\end{equation*}%
est de classe $C^{\infty }$. De plus l'application%
\begin{equation*}
C^{\infty }(M)\longrightarrow C^{\infty }(M^{A},A),f\longmapsto f^{A},
\end{equation*}%
est un homomorphisme injectif d'alg\`{e}bres et on a:%
\begin{eqnarray*}
(f+g)^{A} &=&f^{A}+g^{A} \\
(\lambda \cdot f)^{A} &=&\lambda \cdot f^{A} \\
(f\cdot g)^{A} &=&f^{A}\cdot g^{A}
\end{eqnarray*}%
avec $\lambda \in \mathbb{R}$, $f$ et $g$ appartenant \`{a} $C^{\infty }(M)$.

Lorsque $(a_{\alpha })_{\alpha =1,2,...,\dim (A)}$ est une base de $A$ et
lorsque $(a_{\alpha }^{\ast })_{\alpha =1,2,...,\dim (A)}$ est la base duale
de la base $(a_{\alpha })_{\alpha =1,2,...,\dim (A)}$, l'application%
\begin{equation*}
\sigma :C^{\infty }(M^{A},A)\longrightarrow A\otimes C^{\infty
}(M^{A}),\varphi \longmapsto \dsum\limits_{\alpha =1}^{\dim (A)}a_{\alpha
}\otimes (a_{\alpha }^{\ast }\circ \varphi ),
\end{equation*}%
est un isomorphisme de $A$-alg\`{e}bres. Cet isomorphisme ne d\'{e}pend pas
de la base choisie et l'application%
\begin{equation*}
\gamma :C^{\infty }(M)\longrightarrow A\otimes C^{\infty
}(M^{A}),f\longmapsto \sigma (f^{A}),
\end{equation*}%
est un morphisme d'alg\`{e}bres.

On note $\mathfrak{X}(M^{A})$, l'ensemble des champs de vecteurs sur $M^{A}$%
. Les assertions suivantes sont alors \'{e}quivalentes \cite{b-o}:

\begin{enumerate}
\item $X:C^{\infty}(M^{A})\longrightarrow C^{\infty}(M^{A})$ est un champ de
vecteurs sur $M^{A}$ ;

\item $X:C^{\infty}(M)\longrightarrow C^{\infty}(M^{A},A)$ est une
application lin\'{e}aire v\'{e}rifiant

\begin{equation*}
X(fg)=X(f)\cdot g^{A}+f^{A}\cdot X(g)
\end{equation*}%
pour tous $f$ et $g$ dans $C^{\infty }(M)$.
\end{enumerate}

Ainsi, lorsque%
\begin{equation*}
\theta :C^{\infty }(M)\longrightarrow C^{\infty }(M)
\end{equation*}%
est un champ de vecteurs sur $M$, alors l'application 
\begin{equation*}
\theta ^{A}:C^{\infty }(M)\longrightarrow C^{\infty }(M^{A},A),f\longmapsto 
\left[ \theta (f)\right] ^{A},
\end{equation*}%
est un champ de vecteurs sur $M^{A}$ : le champ de vecteurs $\theta ^{A}$
est le prolongement \`{a} $M^{A}$ du champ de vecteurs $\theta $ sur $M$.

Lorsque $X$ \ est un champ de vecteurs sur $M^{A}$, consid\'{e}r\'{e} comme d%
\'{e}rivation de $C^{\infty }(M)$ dans $C^{\infty }(M^{A},A)$, alors il
existe une d\'{e}rivation et une seule \cite{b-o},%
\begin{equation*}
\widetilde{X}:C^{\infty }(M^{A},A)\longrightarrow C^{\infty }(M^{A},A)
\end{equation*}%
telle que:

\begin{enumerate}
\item $\widetilde{X}$ est $A$-lin\'{e}aire;

\item $\widetilde{X}\left[ C^{\infty}(M^{A})\right] \subset
C^{\infty}(M^{A}) $;

\item $\ \widetilde{X}(f^{A})=X(f)$ pour tout $f\in C^{\infty }(M)$.
\end{enumerate}

L'ensemble $\mathfrak{X}(M^{A})$ des champs de vecteurs sur $M^{A}$\ est
dans ces conditions, un $C^{\infty }(M^{A},A)$-module et une alg\`{e}bre de
Lie sur $A$ \cite{b-o}.

\begin{theoreme}
(de Weil) Si $M$ est une vari\'{e}t\'{e} diff\'{e}rentielle et si $A$ et $B$
sont deux alg\`{e}bres locales, alors l'application%
\begin{equation*}
(M^{A})^{B}\longrightarrow M^{A\otimes B},\eta \longmapsto (id_{A}\otimes
\eta )\circ \gamma
\end{equation*}%
est un isomorphisme de vari\'{e}t\'{e}s diff\'{e}rentielles.
\end{theoreme}

En particulier, on a un isomorphisme entre $TM^{A}$ et $\left( TM\right)
^{A} $.

\bigskip

Pour $x\in M$, $T_{x}M$ d\'{e}signe l'espace tangent en $x$ \`{a} $M$.

On rappelle que la vari\'{e}t\'{e} $M$ est parall\'{e}lisable si son fibr%
\'{e} tangent $TM$ est trivial c'est-\`{a}-dire s'il existe un diff\'{e}%
omorphisme 
\begin{equation*}
\sigma :TM\longrightarrow M\times 
\mathbb{R}
^{n}
\end{equation*}%
tel que le diagramme suivant 
\begin{equation*}
\begin{tabular}{lll}
$TM$ & $\overset{\sigma }{\longrightarrow }$ & $M\times 
\mathbb{R}
^{n}$ \\ 
& \multicolumn{1}{r}{$\underset{\pi _{M}}{\searrow }$} & $\ $\ $\ \
\downarrow pr_{1}$ \\ 
&  & $\ \ \ M$%
\end{tabular}%
\ \ 
\end{equation*}%
commute et que, pour tout $x\in M$ la restriction 
\begin{equation*}
\sigma _{|T_{x}M}:T_{x}M\longrightarrow \left\{ x\right\} \times 
\mathbb{R}
^{n}
\end{equation*}%
soit un isomorphisme d'espaces vectoriels.

Lorsque $(U,\varphi )$ est une carte locale de la vari\'{e}t\'{e} $M$ de
fonctions coordonn\'{e}es $(x_{1},x_{2},...,x_{n})$, l'application 
\begin{equation*}
\ \ \psi :TU^{A}\longrightarrow U^{A}\ \times A^{n},\underset{i=1}{\overset{n%
}{\sum }}\lambda _{i}\cdot \left( \frac{\partial }{\partial x_{i}}\right)
_{\ }^{A}|_{\xi }\ \longmapsto (\xi ,\lambda _{1},...,\lambda _{n})\ 
\end{equation*}%
est un diff\'{e}omorphisme de $A$-vari\'{e}t\'{e}s v\'{e}rifiant $%
pr_{1}\circ \psi =\pi _{M^{A}}$. Ainsi le parall\'{e}lisme local de la vari%
\'{e}t\'{e} $M^{A}$ s'exprime en termes d'existence d'un diff\'{e}omorphisme
de $A$-vari\'{e}t\'{e}s dont la restriction en chaque espace tangent est un
isomorphisme de $A$-modules.

Le but de ce travail est de donner l'\'{e}quivalence du parall\'{e}lisme de $%
M^{A}$ en termes de $A$-vari\'{e}t\'{e}s. On rappelle que lorsque $M$ est
une vari\'{e}t\'{e}, l'alg\`{e}bre de base de $M$ est $C^{\infty }(M)$.
Comme $\mathfrak{X}(M^{A})$ est un $C^{\infty }(M^{A},A)$-module, consid\'{e}%
r\'{e} comme l'ensemble des d\'{e}rivations de $C^{\infty }(M)$ dans $%
C^{\infty }(M^{A},A)$, et est une alg\`{e}bre de Lie sur $A$, et comme $%
M^{A} $ est une $A$-vari\'{e}t\'{e}, ceci signifie que l'alg\`{e}bre de base
de la vari\'{e}t\'{e} $M^{A}$ est $C^{\infty }(M^{A},A)$ et non $C^{\infty
}(M^{A}) $.

\begin{proposition}
La vari\'{e}t\'{e} $M^{A}$ est parall\'{e}lisable si et seulement s'il
existe un diff\'{e}omorphisme de $A$-vari\'{e}t\'{e}s

\begin{equation*}
H:TM^{A}\longrightarrow M^{A}\times A^{n}
\end{equation*}%
tel que le diagramme suivant 
\begin{equation*}
\begin{tabular}{lll}
$TM^{A}$ & $\overset{H}{\longrightarrow }$ & $M^{A}\times A^{n}$ \\ 
& \multicolumn{1}{r}{$\underset{\pi _{M^{A}}}{\searrow }$} & $\ $\ $\ \ \ $\ 
$\ \downarrow pr_{1}$ \\ 
&  & $\ \ $\ $\ \ M^{A}$%
\end{tabular}%
\ 
\end{equation*}%
commute et que pour tout $\xi \in M^{A}$, la restriction%
\begin{equation*}
H|T_{\xi }M^{A}:T_{\xi }M^{A}\longrightarrow \{\xi \}\times A^{n}
\end{equation*}%
soit un isomorphisme de $A$-modules.
\end{proposition}

\begin{demonstration}
$/\Longrightarrow $ Comme la vari\'{e}t\'{e} est parall\'{e}lisable, il
existe un diff\'{e}omorphisme 
\begin{equation*}
\begin{tabular}{lll}
$TM^{A}$ & $\overset{\sigma }{\longrightarrow }$ & $M^{A}\times 
\mathbb{R}
^{n\cdot \dim A}$%
\end{tabular}%
\end{equation*}%
tel que%
\begin{equation*}
\begin{tabular}{lll}
$TM^{A}$ & $\overset{\sigma }{\longrightarrow }$ & $M^{A}\times 
\mathbb{R}
^{n\cdot \dim A}$ \\ 
& \multicolumn{1}{r}{$\underset{\pi _{M^{A}}}{\searrow }$} & $\ $\ $\ \ \ \
\downarrow pr_{1}$ \\ 
&  & $\ \ \ \ \ M^{A}$%
\end{tabular}%
\ 
\end{equation*}%
commute c'est-\`{a}-dire $pr_{1}\circ \sigma =\pi _{M^{A}}$, et que pour
tout $\xi \in M^{A}$, la restriction%
\begin{equation*}
\sigma _{|T_{\xi }M^{A}}:T_{\xi }M^{A}\longrightarrow \{\xi \}\times 
\mathbb{R}
^{n\cdot \dim A}
\end{equation*}%
soit un isomorphisme d'espaces vectoriels.

Soit 
\begin{equation*}
h:A^{n}\longrightarrow 
\mathbb{R}
^{n\cdot \dim A}
\end{equation*}%
un isomorphisme d'espaces vectoriels. Par transport de structure, on munit $%
\mathbb{R}
^{n\cdot \dim A}$ de la structure de $A$-module d\'{e}finie sur $A^{n}$.
Ainsi $h$ devient un isomorphisme de $A$-modules. De la m\^{e}me fa\c{c}on 
\begin{equation*}
\sigma _{|T_{\xi }M^{A}}:T_{\xi }M^{A}\longrightarrow \{\xi \}\times 
\mathbb{R}
^{n\cdot \dim A}
\end{equation*}%
devient un isomorphisme de $A$-modules. En posant $H=(id_{M^{A}}\times
h^{-1})\circ \sigma $, pour tout $\xi \in M^{A}$, on d\'{e}duit que la
restriction 
\begin{equation*}
H|_{T_{\xi }M^{A}}:T_{\xi }M^{A}\longrightarrow \{\xi \}\times A^{n}
\end{equation*}%
est un isomorphisme de $A$-modules.

La diff\'{e}rentiabilit\'{e} de $\sigma $ s'effectuant sur des ouverts de $%
\mathbb{R}
^{2n\cdot \dim A}$, il en est de m\^{e}me pour $H$: ainsi la diff\'{e}%
rentiabilit\'{e} de $H$ s'effectue sur des ouverts de $A^{n}$.

$\Longleftarrow /$ La condition suffisante est \'{e}vidente.
\end{demonstration}

La description du parall\'{e}lisme de $M^{A}$ est la suivante:

\begin{theoreme}
Si $M$ est une vari\'{e}t\'{e} diff\'{e}rentielle de dimension $n$ et si $%
M^{A}$ est la vari\'{e}t\'{e} des points proches de $M$ d'esp\`{e}ce $A$,
alors les assertions suivantes sont \'{e}quivalentes:

\begin{enumerate}
\item La vari\'{e}t\'{e} $M^{A}$ est parall\'{e}lisable;

\item Il existe $n$-champs de vecteurs $X_{1},...,X_{n}$ sur $M^{A}$ tels
qu'en chaque point $\xi \in M^{A}$, $\ $les vecteurs $X_{1}(\xi
),...,X_{n}(\xi )$ forment une base de $T_{\xi }M^{A}$;

\item Le $C^{\infty }(M^{A},A)$-module, $\mathfrak{X}(M^{A})$, des champs de
vecteurs sur $M^{A}$ est un $C^{\infty }(M^{A},A)$-module libre de rang $n$.
\end{enumerate}
\end{theoreme}

\begin{demonstration}
\bigskip Montrons $1/\Longleftrightarrow 2/$

$1/\Longrightarrow 2/$ Comme la vari\'{e}t\'{e} $M^{A}$ est parall\'{e}%
lisable, alors il existe un diff\'{e}omorphisme de $A$-vari\'{e}t\'{e}s

\begin{equation*}
H:TM^{A}\longrightarrow M^{A}\times A^{n}
\end{equation*}%
tel que le diagramme suivant 
\begin{equation*}
\begin{tabular}{lll}
$TM^{A}$ & $\overset{H}{\longrightarrow }$ & $M^{A}\times A^{n}$ \\ 
& \multicolumn{1}{r}{$\underset{\pi _{M^{A}}}{\searrow }$} & $\ $\ $\ \ \
\downarrow pr_{1}$ \\ 
&  & $\ \ $\ $\ \ \ M^{A}$%
\end{tabular}%
\ 
\end{equation*}%
commute et que pour tout $\xi \in M^{A}$, la restriction%
\begin{equation*}
H|_{T_{\xi }M^{A}}:T_{\xi }M^{A}\longrightarrow \{\xi \}\times A^{n}
\end{equation*}%
soit un isomorphisme de $A$-modules.

Pour tout $i=1,2,...,n$, soit $a_{i}=(0,0,...,1_{A},0,...,0)$ o\`{u} $1_{A}$
est \`{a} la $i$-\`{e}me place. Evidemment $(a_{1},a_{2},...,a_{n})$ est une
base du $A$-module $A^{n}$. Pour tout $i=1,2,...,n$, les applications 
\begin{equation*}
\sigma _{i}:M^{A}\longrightarrow \ M^{A}\ \times A^{n},\xi \longmapsto (\xi
,a_{i})
\end{equation*}

et 
\begin{equation*}
X_{i}=H^{-1}\circ \sigma _{i}:M^{A}\longrightarrow \ TM^{A}
\end{equation*}%
sont diff\'{e}rentiables. De plus 
\begin{equation*}
X_{i}:M^{A}\longrightarrow \ TM^{A}
\end{equation*}%
est une section du fibr\'{e} tangent puisque, pour $\xi \in M^{A}$ on a 
\begin{align*}
(\pi _{M^{A}}\circ X_{i})(\xi )& =(pr_{1}\circ H)\left[ (H^{-1}\circ \sigma
_{i})(\xi )\right] \\
& =(pr_{1}\circ H)\left[ H^{-1}(\xi ,a_{i})\right] \\
& =pr_{1}(\xi ,a_{i}) \\
& =\xi \text{.}
\end{align*}%
Ainsi 
\begin{equation*}
\pi _{M^{A}}\circ X_{i}=id_{M^{A}}\text{.}
\end{equation*}%
On conlut que $X_{i}$\ est un champ de vecteurs sur $M^{A}$.

Pour tout $\xi \in M^{A}$, comme $(\xi ,a_{i})_{i=1,2,...,n}$ est une base
du $A$-module $\left\{ \xi \right\} \times A^{n}$, alors $\left[ H^{-1}(\xi
,a_{i})\right] _{i=1,2,...,n}$ est une base du $A$-module $T_{\xi }M^{A}$.
On conclut que les vecteurs $X_{1}(\xi ),...,X_{n}(\xi )$ forment une base
du $A$-module $T_{\xi }M^{A}$.

$2/\Longrightarrow 1/$ On suppose qu'il existe $n$ champs de vecteurs $%
X_{1},...,X_{n}$ sur $M^{A}$ tels qu'en chaque point $\xi \in M^{A}$, $%
(X_{1}(\xi ),...,X_{n}(\xi ))$ soit une base de $T_{\xi }M^{A}$.

L'application 
\begin{equation*}
\ \varphi :M^{A}\ \times A^{n}\ \overset{\ }{\longrightarrow }\ TM^{A},(\xi
,\lambda _{1},...,\lambda _{n})\longmapsto \underset{i=1}{\overset{n}{\sum }}%
\lambda _{i}\cdot X_{i}(\xi )
\end{equation*}%
est un diff\'{e}omorphisme de $A$-vari\'{e}t\'{e}s et 
\begin{equation*}
\ \varphi |_{\{\xi \}\times A^{n}}:\{\xi \}\times A^{n}\ \overset{\ }{%
\longrightarrow }\ T_{\xi }M^{A},(\xi ,\lambda _{1},...,\lambda
_{n})\longmapsto \underset{i=1}{\overset{n}{\sum }}\lambda _{i}X_{i}(\xi )
\end{equation*}%
est un isomorphisme de $A$-modules et sa reciproque%
\begin{equation*}
\ \varphi ^{-1}:\ TM^{A}\longrightarrow M^{A}\ \times A^{n},\underset{i=1}{%
\overset{n}{\sum }}\lambda _{i}X_{i}(\xi )\longmapsto (\xi ,\lambda
_{1},...,\lambda _{n})\ 
\end{equation*}%
est telle que 
\begin{equation*}
pr_{1}\circ \ \varphi ^{-1}=pr_{1}\circ H=\pi _{M^{A}}\text{.}
\end{equation*}%
On conclut alors que la vari\'{e}t\'{e} $M^{A}$ est parall\'{e}lisable.

Montrons $2/\Longleftrightarrow 3/$

$2/\Longrightarrow 3/$ On suppose qu'il existe $n$-champs de vecteurs $%
X_{1},...,X_{n}$ sur $M^{A}$ tels qu'en chaque point $\xi \in M^{A}$, $%
(X_{1}(\xi ),...,X_{n}(\xi ))$ soit une base de $T_{\xi }M^{A}$.

Les champs de vecteurs $X_{1},...,X_{n}$ sont lin\'{e}airement ind\'{e}%
pendants. En effet, si $g_{1},...,g_{n}\in C^{\infty }(M^{A},A)$ sont telles
que%
\begin{equation*}
\underset{i=1}{\overset{n}{\sum }}g_{i}\cdot X_{i}=0\text{,}
\end{equation*}%
alors pour tout $\xi \in M^{A}$, on a%
\begin{equation*}
\underset{i=1}{\overset{n}{\sum }}g_{i}(\xi )\cdot X_{i}(\xi )=0\text{.}
\end{equation*}%
Comme $(X_{1}(\xi ),...,X_{n}(\xi ))$ est une base de $T_{\xi }M^{A}$, ainsi 
$g_{i}(\xi )=0$ pour tout $i=1,2,...,n$. Comme $\xi $ est quelconque, on
conlut que $g_{i}=0$ pour tout $i=1,2,...,n$.

La famille $X_{1},...,X_{n}$ engendre $\mathfrak{X}(M^{A})$, en effet, si $%
Y\in \mathfrak{X}(M^{A})$ et $\xi \in M^{A}$, on a :%
\begin{equation*}
Y\left( \xi \right) =\underset{i=1}{\overset{n}{\sum }}\lambda _{i}\cdot
X_{i}(\xi )
\end{equation*}%
avec les $\lambda _{i}\in A$.

L'application $\ $ 
\begin{equation*}
\ M^{A}\overset{Y}{\longrightarrow }\ TM^{A}\ \overset{H}{\longrightarrow }%
M^{A}\times A^{n}\ \overset{pr_{2}}{\longrightarrow }\ A^{n}\overset{pr_{i}}{%
\longrightarrow }A,\xi \longmapsto \lambda _{i},
\end{equation*}%
est diff\'{e}rentiable. En posant $f_{i}=pr_{i}\circ pr_{2}\circ H\circ Y$,
on a $f_{i}(\xi )=\ \lambda _{i}$ et%
\begin{align*}
Y\left( \xi \right) & =\underset{i=1}{\overset{n}{\sum }}f_{i}(\xi )\cdot
X_{i}(\xi ) \\
& =\ \left( \underset{i=1}{\overset{n}{\sum }}f_{i}\cdot X_{i}\right) (\xi )%
\text{.}
\end{align*}%
Comme $\xi $ est quelconque, alors 
\begin{equation*}
Y=\underset{i=1}{\overset{n}{\sum }}f_{i}\cdot X_{i}
\end{equation*}

Ainsi $X_{1},...,X_{n}$ une base du $C^{\infty }(M^{A},A)$-module $\mathfrak{%
X}(M^{A})$. On conlut que $\mathfrak{X}(M^{A})$ est un $C^{\infty }(M^{A},A)$%
-module libre de rang $n$.

$3/\Longrightarrow 2/$ On suppose que $\mathfrak{X}(M^{A})$ est un $%
C^{\infty }(M^{A},A)$-module libre de rang $n$. Soit $(X_{1},...,X_{n})$ une
base du $C^{\infty }(M^{A},A)$-module $\mathfrak{X}(M^{A})$.
\end{demonstration}

Si $\alpha _{1}(\xi ),...,\alpha _{n}(\xi )$ sont des \'{e}l\'{e}ments de $A$
tels que 
\begin{equation*}
\underset{i=1}{\overset{n}{\sum }}\alpha _{i}(\xi )\cdot X_{i}\left( \xi
\right) =0
\end{equation*}%
pour tout $\xi \in M^{A}$, pour tout $i=1,2,...,n$, soit 
\begin{equation*}
f_{i}:M^{A}\longrightarrow A,\xi \longmapsto \alpha _{i}(\xi )\text{.}
\end{equation*}%
Pour $\eta \in M^{A}$, il existe $Y\in \mathfrak{X}(M^{A})$ tel que $Y(\eta
)=\underset{i=1}{\overset{n}{\sum }}f_{i}(\eta )\cdot X_{i}\left( \eta
\right) $. Comme $Y$ est diff\'{e}rentiable au voisinage de $\eta $, il en
est de m\^{e}me de $f_{i}$ au voisinage de $\eta $. Comme $\eta $ est
quelconque, on d\'{e}duit que les $f_{i}$ sont diff\'{e}rentiables. Ainsi,
on a 
\begin{eqnarray*}
0 &=&\underset{i=1}{\overset{n}{\sum }}\alpha _{i}(\xi )\cdot X_{i}\left(
\xi \right) \\
&=&\underset{i=1}{\overset{n}{\sum }}f_{i}(\xi )\cdot X_{i}\left( \xi \right)
\\
&=&\left( \underset{i=1}{\overset{n}{\sum }}f_{i}\cdot X_{i}\right) \left(
\xi \right)
\end{eqnarray*}%
pour tout $\xi \in M^{A}$. Comme les champs de vecteurs $X_{1},...,X_{n}$
forment une base du $C^{\infty }(M^{A},A)$ module $\mathfrak{X}(M^{A})$,
alors 
\begin{equation*}
f_{1}=...=f_{n}=0\text{.}
\end{equation*}%
C'est-\`{a}-dire que pour tout $\xi \in M^{A}$, les $\alpha _{i}(\xi )=0$.
On conclut que la famille $\left( X_{1}(\xi ),...,X_{n}(\xi )\right) $ est
libre pour tout $\xi \in M^{A}$.

\begin{demonstration}
De plus, pour $v\in T_{\xi }M^{A}$, il existe un champ de vecteurs $Y\in 
\mathfrak{X}(M^{A})$ telle que $Y(\xi )=v.$ Puisque 
\begin{equation*}
Y=\underset{i=1}{\overset{n}{\sum }}f_{i}\cdot X_{i}
\end{equation*}%
$\ \ $o\`{u} chaque $f_{i}\in C^{\infty }(M^{A},A)$, alors, 
\begin{equation*}
v=\underset{i=1}{\overset{n}{\sum }}f_{i}(\xi )\cdot X_{i}(\xi )\text{.}\ \ 
\end{equation*}%
Ainsi, la famille $\left( X_{1}(\xi ),...,X_{n}(\xi )\right) $ engendre le $%
A $-module $T_{\xi }M^{A}$.

On conclut alors qu'en chaque point $\xi \in M^{A}$, les vecteurs $X_{1}(\xi
),...,X_{n}(\xi )$ forment une base du $A$-module $T_{\xi }M^{A}$.
\end{demonstration}

\begin{corollaire}
Si $M$ est une vari\'{e}t\'{e} parall\'{e}lisable, alors la vari\'{e}t\'{e}
des points proches $M^{A}$ est parall\'{e}lisable.
\end{corollaire}

\end{document}